\documentclass[a4paper,twoside]{amsart}
\input amssym.def
\input amssym.tex
\font\tencmmib=cmmib10
\skewchar\tencmmib='177
\font\sevencmmib=cmmib7
\skewchar\sevencmmib='177
\font\fivecmmib=cmmib5
\skewchar\fivecmmib='177
\newfam\cmmibfam
\textfont\cmmibfam=\tencmmib
\scriptfont\cmmibfam=\sevencmmib
\scriptscriptfont\cmmibfam=\fivecmmib
\font\tencmbsy=cmbsy10
\skewchar\tencmbsy='60
\font\sevencmbsy=cmbsy7
\skewchar\sevencmbsy='60
\font\fivecmbsy=cmbsy5
\skewchar\fivecmbsy='60
\newfam\cmbsyfam
\textfont\cmbsyfam=\tencmbsy
\scriptfont\cmbsyfam=\sevencmbsy
\scriptscriptfont\cmbsyfam=\fivecmbsy 
\newenvironment{dem}{\textit{Proof :}}{$\Box$\newline}

\newtheorem{Theorem}{Theorem}
\newtheorem{Lemma}{Lemma}

\newtheorem{Remark}{Remark}

\def\R{{\mathbb{R}}} 
\def\N{{\mathbb{N}}}

\def\H{{\mathbb{H}}}
\def\epsilon{{\varepsilon}}
\def\phi{{\varphi}}
\def\theta{{\vartheta}}

\def\cal{{}}

\DeclareMathOperator{\graph}{graph}
\DeclareMathOperator{\argth}{Argth}

\def\ds{\displaystyle}
\begin{document}

\title{Entire spacelike radial graphs in the Minkowski space,
asymptotic to the light-cone, with prescribed scalar curvature}
\author{Pierre Bayard \& Philippe Delano\" e}
\thanks{The first author was supported by the project UNAM-PAPITT IN 101507; the second author is supported by the CNRS.}
\address{Pierre Bayard, Instituto de
F\'{\i}sica y Matem\'aticas, U.M.S.N.H. Ciudad Universitaria, CP.
58040 Morelia, Michoac\'an, Mexico}
\email{bayard@ifm.umich.mx}
\address{Philippe Delano\" e, Universit\'e de
Nice--Sophia Antipolis, Laboratoire J.--A. Dieudonn\'e, Parc Valrose, 06108 Nice Cedex 2, France}
\email{Philippe.DELANOE@unice.fr}
\date{} 
\begin{abstract}
Existence and uniqueness in ${\Bbb R}^{n,1}$ of entire spacelike
hypersurfaces contained in the future of the origin $O$ and asymptotic
to the light-cone,
with scalar curvature prescribed at their generic point $M$ as a
negative function of the unit vector
$\overrightarrow{Om}$ pointing in the direction of
$\overrightarrow{OM}$, divided by the square of the norm of
$\overrightarrow{OM}$ (a dilation invariant problem). The solutions are seeked as graphs
over the future unit-hyperboloid emanating from $O$ (the hyperbolic space); radial upper and lower solutions are constructed which, relying on a previous
result in the Cartesian setting, imply their existence.\\

\noindent R{\tiny ESUME}. Existence et unicit\'e dans ${\Bbb R}^{n,1}$ d'hypersurfaces enti\`eres de genre espace contenues dans le futur de l'origine $O$
et asymptotes au c\^one de lumi\`ere, dont la courbure scalaire est prescrite au point g\'en\'erique $M$  comme fonction n\'egative du vecteur unit\'e
$\overrightarrow{Om}$ pointant en direction de $\overrightarrow{OM}$, divis\'ee par le carr\'e de la norme du vecteur $\overrightarrow{OM}$ (un probl\`eme
invariant par homoth\'etie). Les solutions sont cherch\'ees comme graphes sur l'hyperbolo\"{\i}de-unit\'e futur \'emanant de $O$ (l'espace hyperbolique);
des solutions sup\'erieure et inf\'erieure radiales sont construites qui, d'apr\`es un r\'esultat ant\'erieur en cart\'esien, impliquent l'existence de
telles solutions.
\end{abstract}

\maketitle
\thispagestyle{empty}
\markboth{\textsl{P. Bayard \& Ph. Delano\"e}}
{\textsl{Entire spacelike radial graphs with prescribed scalar curvature}}
\section*{Introduction} 

The Minkowski space $\R^{n,1}$ is the affine Lorentzian manifold $\R^n\times\R$ endowed with the metric
$$ds^2=d{X'}^2-dX_{n+1}^2\ ,\ {\rm where}\ d{X'}^2=dX_1^2+\ldots+dX_n^2\ ,$$ 
setting $X=(X',X_{n+1})\in \R^n\times\R$, and time-oriented by $dX_{n+1}>0$.
Distinguishing the origin $O$ of $\R^{n,1},$ let
$$\H=\{x\in\R^{n,1}|\ \vert \overrightarrow{Ox} \vert^2=|x'|^2-{x_{n+1}}^2=-1,\ x_{n+1}>0\}\ ,$$ 
be the future unit-hyperboloid, model of the hyperbolic space in $\R^{n,1}.$
If $\varphi$ is a real function defined on $\H,$ we define the \emph{radial
graph} of $\varphi$ by
$$\graph_{\H}\varphi=\{X\in \R^{n,1}, \overrightarrow{OX}=e^{\varphi(x)} \overrightarrow{Ox},\ x\in\H\}\ .$$
This is a hypersurface contained in the future open solid cone
$$C^+=\{X\in\R^{n,1}|\ X_{n+1}>|X'|\}\ .$$
We say that $\varphi$ is spacelike if its graph is a spacelike
hypersurface, which means that the metric induced on it is Riemannian. Conversely, a spacelike and connected hypersurface in $C^+$ is the radial graph of a
uniquely determined function $\varphi:\H\rightarrow\R.$ Of course, such a graph may also be considered as the Cartesian graph of some function
$u:\R^n\rightarrow\R$ 
$$\graph_{\R^n}u=\{(x',u(x')),\ x'\in\ \R^n\},$$
and the correspondence between the two representations is bijective passing from the Cartesian chart $X=(X',X_{n+1})$ restricted to $C^+$, to the polar
chart
$(x,\rho)\in \H\times (0,\infty)$ of $C^+$ defined by:
$$\rho=\sqrt{-\vert\overrightarrow{OX}\vert^2},\ \overrightarrow{Ox}=\frac{1}{\rho}\overrightarrow{OX}\ .$$
Recall that the principal curvatures $(\kappa_1,\ldots,\kappa_n)$ at a point of a spacelike hypersurface are the eigenvalues of its curvature endomorphism
$dN,$ where
$N$ is the future oriented unit normal field, and the $m^{th}$ mean curvature (denoted by $H_m$) is the $m^{th}$ elementary symmetric function of its
principal curvatures: $H_m=\sigma_m(\kappa_1,\ldots,\kappa_n).$
For each real $\lambda>0$, the cone $C^+$ is globally invariant under the ambient dilation $X\mapsto \lambda X$ of $\R^{n,1}$ and the above $m$-th mean
curvature is $(-m)$-homogeneous; specifically, it transforms like $H_m(\lambda X)=\lambda^{-m}H_m(X).$ It is thus natural to pose, as in \cite[Theorem
1]{del1}, the following inverse problem for $H_m$: given a positive function $h>0$ on $\H$ tending to 1 at infinity, find a spacelike hypersurface
$\Sigma$ in
$C^+$, asymptotic to $\partial C^+$ at infinity, such that, for each point $X\in \Sigma$, the $m$-th mean curvature of $\Sigma$ at $X$ is given by:
\begin{equation}\label{eqn k curv hom}
\frac{1}{{n \choose
m}}H_m(X)=\frac{1}{(-|\overrightarrow{OX}|^2)^{\frac{m}{2}}}[h(x)]^m\ , {\rm with}\
\overrightarrow{Ox}=\frac{\overrightarrow{OX}}{\sqrt{-|\overrightarrow{OX}|^2}}\ .
\end{equation}
By construction, this problem is dilation invariant; moreover, as explained below, the positivity of $h$ makes it elliptic. Actually, introducing the positivity cone
\cite{ivo1} of $\sigma_m$:
$$\Gamma_m=\{\kappa\in \R^n,\ \forall i=1,\ldots,m,\ \sigma_i(\kappa)>0\}\ ,$$
and recalling McLaurin's inequalities (satisfied on $\Gamma_m$):
$$0<(H_m)^{\frac{1}{m}} \leq (H_{m-1})^{\frac{1}{m-1}} \leq \ldots \leq H_2^{\frac{1}{2}}\leq H_1\ ,$$
we note that, if a hypersurface $\Sigma=\graph_{\R^n}u$ solves 
(\ref{eqn k curv hom}) with the asymptotic condition, then the time-function $u$ must assume a minimum on $\Sigma$  and, as readily checked (using
\textit{e.g.}
\cite[p.245]{bay2}), the principal curvatures of $\Sigma$ at such a minimum point of $u$ must lie in $\Gamma_m$. Now equation (\ref{eqn k curv hom})
combined with McLaurin's inequalities forces the principal curvatures of
$\Sigma$ to stay in
$\Gamma_m$ \textit{everywhere}. Let us call any spacelike hypersurface of
$C^+$ having this property,
$m$-admissible; accordingly, a function $\phi:\H
\to \R$ (resp.
$u:\R^n \to
\R$) is called $m$-admissible, provided $\graph_{\H}\phi$  (resp. $\graph_{\R^n}u$) is so. The condition of $m$-admissibility is local (and open); one may
thus speak of a function $\phi:\H \to \R$ being $m$-admissible \textit{at a point} (hence nearby) whenever $\graph_{\H}\phi$ is so at that point. We will
seek the solution hypersurface $\Sigma$ as the radial graph of some
$m$-admissible function
$\varphi:\H \to \R$ vanishing at infinity (to comply with the asymptotic condition). Equation (\ref{eqn k curv hom}) then reads
\begin{equation}
\label{eqfmphi}
F_m(\varphi)=h,
\end{equation}
with the radial operator $F_m$ defined by:
$$F_m(\varphi)=e^{\phi}\left[\frac{1}{{n \choose m}}H_m(X)\right]^{\frac{1}{m}},\ X\in\graph_{\H}\phi.$$
For briefness, we will not compute here explicitely the general expression of the operator $F_m$ (keeping it for a further study) -- its restriction to
radial functions will suffice (see section 3.3 below). We will rely instead on
the well-known corresponding Cartesian expression (see
\textit{e.g.}
\cite{bay1}) combined with a few basic properties of
$F_m$ recorded in the next section (and proved with elementary arguments).\\
Furthermore, we will essentially restrict to the case $m=2$ (and freely say 'admissible', for short, instead of
'2-admissible'). Since 
$H_2$ is related to the scalar curvature $S$ by
$S=-2H_2$,
our present study is really about the prescription of the scalar curvature, at a generic point $X$ of a radial graph, as a negative function of $x\in \H$
(with $x$ given as in (\ref{eqn k curv hom})) divided by the square of the norm of
$\overrightarrow{OX}.$ Aside from the origin $O$ of the ambient space $\R^{n,1}$, we will distinguish a point $o$ in $\H$ and set $r=r(x)$ for the
hyperbolic distance from $o$ to $x\in \H$; accordingly, a function on $\H$ will be called \textit{radial} whenever it factors through a function of $r$
only. Our main result is the following:
\begin{Theorem}\label{mainth}
For $\alpha\in (0,1)$, let $h:\H\rightarrow (0,\infty)$ be a function of class $C^{2,\alpha}$ with
$$\lim_{r(x)\rightarrow+\infty} h(x)=1\ .$$
Assume that the functions $h^-$ and $h^+$ defined on
$\R^+$ by
$$h^-(r)=\sup_{r(x)=r} h(x)\mbox{ and }h^+(r)=\inf_{r(x)=r} h(x)$$
satisfy 
$$ \int_0^{+\infty}(h^--1)_+dr<+\infty\ ,\ \int_0^{+\infty}(1-h^+)_+dr<+\infty\ ,$$  
where $(h^--1)_+$ (resp. $(1-h^+)_+$) means the positive part of $h^--1$ (resp. $1-h^+$).
Then the equation 
$$F_2(\varphi)=h$$
has a unique admissible solution of class $C^{4,\alpha}$ such that $\lim_{r(x)\rightarrow+\infty}\varphi(x)=0.$ 
\end{Theorem}
\begin{Remark}
{\rm From Lemma \ref{limite phi finie} below,  anytime the function $h$ is radial, the integral convergence conditions of Theorem \ref{mainth} appears
necessary for the existence of bounded solutions.}
\end{Remark}
An analogous problem in the Euclidean setting is solved for the Gauss curvature in \cite[Th\'eor\`eme 1]{del1}, and in \cite{tw,cns4} some related problems
are studied. In the Lorentzian setting, the prescription of the mean curvature for entire graphs is studied in \cite{bar} and that of the Gauss curvature in
\cite{li,gujisch,basch}. In
\cite{bay2}, the scalar curvature is prescribed in Cartesian coordinates $x_{n+1}=u(x_1,\ldots,x_n).$
\\
The outline of the paper is as follows. In section 1, we prove that there exists at most one solution vanishing at infinity for equation
(\ref{eqfmphi}) with $m\in \{1,\ldots,n\}$. In section 2, relying on \cite{bay2}, we prove the existence of a solution when $m=2$ , provided upper and lower
barriers are known. The latter are constructed, as radial functions, in section 3.
\section{Uniqueness}\label{uniqueness}
We first require a few basic properties of the operator $F_m$. It is a nonlinear second order scalar differential operator defined on
$m$-admissible real functions on
$\H$. The dilation invariance of (\ref{eqn k curv hom}) implies the identity:
\begin{equation}
\label{uptoc}
F_m(\psi+c)\equiv F_m(\psi)\ ,
\end{equation}
for every $m$-admissible function $\psi:\H\rightarrow\R$ and constant $c$; linearizing at $\psi$ yields
$$dF_m(\psi)(1)\equiv 0\ .$$
Furthermore, we have:
\begin{Lemma}\label{elliplem}
For each $m$-admissible function $\psi$, the linear differential operator $dF_m(\psi)$ is elliptic everywhere on $\H$, with positive-definite symbol.
\end{Lemma}
Summarizing for later use, the expression of $dF_m(\psi)$, in the chart $x'\in {\Bbb R}^n$ of ${\Bbb H}$, at a fixed $m$-admissible function $\psi$ reads
like:
\begin{equation}\label{expresslinear}
\delta \psi \mapsto dF_m(\psi)(\delta \psi)= \sum_{1\leq i,j\leq n} B_{ij} \frac{\partial^2}{\partial {x'}_i\partial {x'}_j}(\delta \psi)+ \sum_{i=1}^n B_i
\frac{\partial}{\partial {x'}_i}(\delta \psi)\ ,
\end{equation}
with the $n\times n$ matrix $(B_{ij})$ symmetric positive definite (depending on $\psi$, of course, like the $B_i$'s). We now proceed to proving Lemma
\ref{elliplem}.\\
\begin{dem}
We require the Cartesian operator $v\mapsto G_m(v):=F_m(\psi)$ defined on
$m$-admissible functions $v:\R^n \to \R$ by:
\begin{equation}\label{samegraph}
\graph_{\R^n}v=\graph_{\H}\psi\ .
\end{equation}
The ellipticity of $dG_m(v)$ and the positive-definiteness of its symbol are well-known \cite{ivo2,tru,bay1}. Its expression thus starts out like
$$dG_m(v)(\delta v)= \sum_{1\leq i,j \leq n} A_{ij} \frac{\partial^2}{\partial {X'}_i\partial {X'}_j}(\delta v)+\ \mbox{{\rm lower order terms}}\ ,$$
with the matrix $(A_{ij})$ symmetric positive  definite. The $m$-admissible function $\psi$ on
$\H$ such that (\ref{samegraph}) holds, is related to $v$, in the chart $x'=(x_1,\ldots,x_n)\in \R^n$, by:
$$v(X')=\sqrt{1+\vert x' \vert^2} \exp{[\psi(x')]},\ {\rm with}\ \overrightarrow{OX'}= e^{\psi(x')}\overrightarrow{Ox'}\ .$$
Varying $\psi$ by $\delta \psi$ thus yields for the corresponding variation $\delta v$ of $v$ the following expression: $\delta v(X')=w(X') \delta \psi
(x')$, with $\ds w(X')=\left[v-\sum_{i=1}^n X'_i\frac{\partial v}{\partial X'_i}\right](X')$. Since the graph lies in $C^+$ and it is spacelike, we
have $v(X')>\vert X' \vert$ and (using Schwarz inequality) $\sum_{i=1}^n X'_i\frac{\partial v}{\partial X'_i}<\vert X' \vert$, therefore $w>0$. Moreover,
up to lower order terms, we have:
$$ \frac{\partial^2}{\partial {X'}_i\partial {X'}_j}(\delta v)(X')=w(X')\sum_{1 \leq i,j \leq n}  \frac{\partial^2}{\partial {x'}_k\partial {x'}_l}(\delta
\psi)(x') \frac{\partial x'_k}{\partial X'_i} \frac{\partial x'_l}{\partial X'_j}$$
with $\ds x'_k=\frac{X'_k}{\sqrt{v^2(X')-\vert X' \vert^2}}$. We thus find in (\ref{expresslinear}):
$\ds B_{kl}=w(X') \sum_{1\leq i,j \leq n} A_{ij}  \frac{\partial x'_k}{\partial X'_i} \frac{\partial x'_l}{\partial X'_j}$
and the ellipticity of $\delta \psi\mapsto
dF_m(\psi)(\delta
\psi)$ follows.
\end{dem}

We need also a more specific (ellipticity) property of the operator $F_m$, namely:
\begin{Lemma}
\label{admisegment}
For each couple $(\varphi_0,\varphi_1)$ of $m$-admissible real functions on $\H$ and each point $x_0\in \H$ where $\varphi=\varphi_1-\varphi_0$ assumes a
local extremum, the whole segment $t\in [0,1]\to \varphi_t=\varphi_0+t\varphi$ consists of $m$-admissible functions at the point $x_0$.
\end{Lemma}
\begin{dem}
The analogue of Lemma \ref{admisegment} is fairly standard in the Cartesian setting, using the expression of the operator $G_m$ introduced
in the proof of Lemma \ref{elliplem} (see \cite{bay1}) together with the well-known fact: $\forall \kappa \in \Gamma_m, \forall
i\in \{1,\ldots,n\},\ \frac{\partial \sigma_m}{\partial \kappa_i}(\kappa) >0$. Here, we will simply reduce the proof to that setting (and let the reader
complete the argument). Let us first normalize the situation at an extremum point
$x_0\in
\H$ of
$\varphi$. From (\ref{uptoc}), we may assume
$\varphi(x_0)=0$. Moreover, we may assume that $\varphi$ has a local minimum at $x_0$ (if not, switch $\varphi_0$ and $\varphi_1$). Finally, setting
$\graph_{\H}\varphi_a=\graph_{\R^n}u_a$ for $a=0,1$, and performing if necessary a  suitable Lorentz transform (hyperbolic rotation), we may take
$x_0=(0,1)\in \R^n\times\R$ thus with $u_a(0)=1$. For $t\in [0,1]$ and near $x_0$, set $\Sigma_t=\graph_{\R^n}u_t$ for the hypersurface
$\graph_{\H}\varphi_t$. We must prove that $\Sigma_t$ is $m$-admissible at $x_0$. For $X_t\in \R^{n,1}$ lying in $\Sigma_t$, we have:
$\overrightarrow{OX_t}=e^{t\varphi(x)}\overrightarrow{OX_0}$ with $\ds \overrightarrow{Ox}=\frac{\overrightarrow{OX_0}}{\sqrt{-\vert \overrightarrow{OX_0}
\vert^2}}$. In the Cartesian setting, we thus have (sticking to the $\R^n$-valued charts used in the preceding proof):
$$u_t(X_t')=e^{t\varphi(x')}u_0[e^{-t\varphi(x')}X_t']\ ,$$
here with $\ds x'=\frac{X'_0}{\sqrt{u_0^2(X'_0)-\vert X'_0\vert^2}},\ X_t'=e^{t\varphi(x')}X'_0$, and $(X'_0,u_0(X'_0))\in \graph_{\R^n}u_0$;
moreover, the lemma boils down to proving that $u_t$ is $m$-admissible at $X'_t=0$. A routine calculation yields at $X'_t=0$ the equalities:
$$\frac{\partial u_t}{\partial X'_{ti}}(0)=\frac{\partial u_0}{\partial X'_{0i}}(0),\ \frac{\partial^2 u_t}{\partial X'_{ti}\partial
X'_{tj}}(0)=\frac{\partial^2 u_0}{\partial X'_{0i}
\partial X'_{0j}}(0)+ t
\frac{\partial^2 \varphi}{\partial x'_i\partial x'_j}(0)\ ,$$
where, in the second one, the matrix $\ds \left[\frac{\partial^2 \varphi}{\partial x'_i\partial x'_j}(0)\right]_{1\leq i,j\leq n}$ is non-negative. The rest
of the proof is now standard, thus omitted.
\end{dem}
\begin{Theorem}
\label{uniqth}
The operator $F_m$ is one-to-one on $m$-admissible functions of class $C^2$ vanishing at infinity.
\end{Theorem}
\begin{dem}
Let us argue by contradiction. Let $\varphi_0,\varphi_1$ be two $m$-admissible $C^2$ functions vanishing at infinity and having the same image by $F_m$. For
$t\in [0,1]$, set
$\varphi_t=\varphi_0+t\varphi$ with
$\varphi=\varphi_1 - \varphi_0$. Since $\varphi$ vanishes at infinity, if $\varphi \not\equiv 0$, it assumes a nonzero local extremum (a maximum, say, with
no loss of generality) at some point $x_0\in \H$. By Lemma \ref{admisegment}, the whole segment $t\in[0,1]\to \varphi_t$ is $m$-admissible in a
neighborhood
$\Omega$ of
$x_0$ where
$\varphi$ thus satisfies the second order linear equation $L\varphi=0$ with $L1=0$ and the operator $L$ given by $L=\int_0^1dF_m(\varphi_t) dt$. Combining
Lemma
\ref{elliplem} above with Hopf's strong Maximum Principle (see \cite{gt}), we get $\varphi\equiv \varphi(x_0)$ throughout $\Omega$. By
connectedness, we infer
$\varphi\equiv \varphi(x_0)\not= 0$ on the whole of $\H$, contradicting $\lim_{r(x)\rightarrow+\infty}\varphi = 0$. So, indeed, we must have $\varphi
\equiv 0$, in other words $F_m$ is one-to-one.
\end{dem}

\section{Existence of a solution reduced to that of upper and lower solutions}
\begin{Theorem}
Let $h:\H\rightarrow\R$ be a function of class $C^{2,\alpha},$ for some $\alpha \in (0,1)$, such that  there exists $\varphi^-\in C^{4,\alpha}(\H)$ with
$\graph_{\H}\phi^-$ strictly convex and spacelike, and $\varphi^+\in C^{2}(\H)$ with $\graph_{\H}\phi^+$ spacelike, satisfying
$$F_2(\varphi^-)\geq h,\ F_2(\varphi^+)\leq h\mbox{ and }\lim_{r(x)\rightarrow+\infty}\varphi^{\pm}=0.$$
Then the equation 
$$F_2(\varphi)=h$$
has a unique admissible solution of class $C^{4,\alpha}$ such that $\lim_{r(x)\rightarrow+\infty}\varphi(x)=0.$ Moreover $\varphi$ satisfies the pinching:
$$\varphi^-\leq \varphi\leq \varphi^+.$$
\end{Theorem}
\begin{Remark}
{\rm Since $\varphi$ is a bounded function,  the hypersurface $M=\graph_{\H}(\varphi)$ is entire. More precisely, denoting by $\varphi_{min}$ and
$\varphi_{max}$ two constants such that
$\varphi_{min}\leq\varphi\leq\varphi_{max},$ the
function $u:\R^n\rightarrow\R$ such that $\graph_{\R^n}(u)=\graph_{\H}(\varphi)$ satisfies $u_{min}\leq
u\leq u_{max}$ where $u_{min}$ (resp. $u_{max}$) is such that  $\graph_{\R^n}(u_{min})=\graph_{\H}(\varphi_{min})$ (resp.
$\graph_{\R^n}(u_{max})=\graph_{\H}(\varphi_{max})).$ Noting that the graphs of $u_{min}$ and $u_{max}$ are hyperboloids,  we see that the inequality
$u\geq u_{min}$ implies that $M$ is entire, and the inequality $u\leq u_{max}$ implies that $M$ is asymptotic to the lightcone.}
\end{Remark}
\begin{dem}
The asserted uniqueness follows from Theorem \ref{uniqth}; so let us focus on the existence part. A straightforward comparison principle, using
(\ref{expresslinear}) and Lemma
\ref{admisegment}, implies
$\varphi^-\leq
\varphi^+$ on
$\H.$ Let
$u^-,\ u^+:\R^n\rightarrow\R$ be such that $\graph_{\R^n}(u^\pm)=\graph_{\H}(\varphi^\pm).$ Set $H$ for the function on $\R^{n,1}$ defined by: 
\begin{equation}\label{Definition H}
H(X)=\frac{(^n_2)}{|X_{n+1}|^2-|X'|^2}\left[h\left(\frac{X}{\sqrt{|X_{n+1}|^2-|X'|^2}}\right)\right]^2.
\end{equation}
The spacelike functions $u^-$ and $u^+$ satisfy:
$$H_2[u^-]\geq H(.,u^-),\ H_2[u^+]\leq H(.,u^+),\ u^-\leq u^+\mbox{ and }\lim_{|x'|\rightarrow \infty}[u^\pm(x')-|x'|]= 0\ ,$$
where $H_2[u^{\pm}]$ stands for the second mean curvature of the graph of $u^{\pm}.$
Theorem 1.1 in \cite{bay2} asserts the existence of a function $u:\R^n\rightarrow\R,$ belonging to $C^{4,\alpha},$ spacelike, such that $H_2[u]=H(.,u)$ in
$\R^n,$ $\displaystyle{\lim_{|x'|\rightarrow+\infty}u(x')-|x'|=0,}$ and $u^-\leq u\leq u^+.$  The function $\varphi:\H\rightarrow\R$ such that
$\graph_{\H}(\varphi)=\graph_{\R^n}(u)$ is a solution of our original problem. 
\end{dem}

\section{Construction of radial upper and lower solutions}
In the sequel of the paper, we first solve the Dirichlet problem on a bounded set in $\H$ (section \ref{parag 1}) then proceed to proving the existence and
uniqueness of an entire solution in the radial case and study its properties (sections \ref{parag 2} and \ref{parag 3}); finally, we construct the required
radial barriers (section
\ref{parag 4}).
\subsection{The Dirichlet problem}\label{parag 1}
\begin{Theorem} \label{theorem Dirichlet problem}
Given $\alpha\in (0,1)$, let $\Omega$ be a uniformly convex bounded open subset of $\H$ with $C^{2,\alpha}$ boundary, $h:\Omega\rightarrow\R$ be a positive
function of class
$C^{2,\alpha},$ and $\varphi_0:\overline{\Omega}\rightarrow\R$ be a spacelike function of class $C^{2,\alpha}$ whose radial graph is strictly convex. Then
the Dirichlet problem
\begin{equation}\label{Dirichlet problem}
F_2(\varphi)=h\mbox{ in }\Omega,\ \varphi=\varphi_0\mbox{ on }\partial\Omega,
\end{equation}
has a unique admissible solution of class $C^{4,\alpha}.$
\end{Theorem}
\begin{dem}
The proof of uniqueness follows the lines of the proof of Theorem \ref{uniqth}; let us focus on
the existence part. Setting $x=(x',\sqrt{1+\vert x' \vert^2})\in\R^n\times\R,$ and
$$\Omega'=\{e^{\varphi_0(x)}x',\ x\in\Omega\},\ 
u_0(e^{\varphi_0(x)}x')=e^{\varphi_0(x)}\sqrt{1+\vert x' \vert^2}\ ,$$
problem (\ref{Dirichlet problem}) is equivalent to the Dirichlet problem: 
\begin{equation}
H_2[u]=H(.,u)\mbox{ in }\Omega',\ u=u_0\mbox{ on }\partial\Omega',
\end{equation}
where $H_2$ is the scalar curvature operator acting on spacelike graphs defined on $\Omega'\subset \R^n,$ and $H$ is defined on $\Omega'\times\R$ by
(\ref{Definition H}). We know essentially from 
\cite{bay1,ur} that this problem is solvable (with an adaptation here because the function $H$ depends also on $u;$ the existence is proved by a
classical fixed point argument \cite{gt} and the required \textit{a priori} estimates are carried out in \cite[p.251]{bay2}).
\end{dem}
\subsection{Existence and uniqueness of entire radial solutions}\label{parag 2}
The aim of this section is to prove the following result :
\begin{Theorem}\label{theorem entire solution}
For $\alpha\in (0,1)$, let $h:\R^+\rightarrow\R$ be a positive function of class $C^{2,\alpha}$ constant on some neighborhood of 0 and let
$\varphi_0$ be a real number.
Recall $r=r(x)$ denotes the hyperbolic distance of $x\in \H$ from a fixed origin $o\in \H$. The problem:
\begin{equation}\label{eqn k phi}
F_2(\varphi)(x)=h(r) \mbox{ for all }x\in\H,\ \varphi(o)=\varphi_0,
\end{equation}
admits a unique admissible radial solution $\varphi:\H\rightarrow\R$ of class $C^{4,\alpha}.$
\end{Theorem}
\begin{dem}
\textit{Existence:} let $B_i$ denote the ball in $\H$ with center $o$ and radius $i\in \N^*$, and $\varphi_{i}$ be the admissible solution of the Dirichlet
problem:
\begin{equation}\label{radial Dirichlet problem i}
F_2(\varphi)=h,\ \varphi_{|\partial B_i}=0,
\end{equation}
given by Theorem \ref{theorem Dirichlet problem}. By radial symmetry and uniqueness, $\varphi_i$ is a radial function: $\varphi_i(x)=f_i(r)$ for some
function
$f_i:[0,i]\rightarrow\R$. By uniqueness again, for $j>i$, the function $\varphi_{j}-\varphi_{i}$ must be constant on $B_i$. Therefore
$f_{j}'(r)\equiv f_{i}'(r)$ for
$r\in [0,i].$ We may thus define $g$ on $\R^+$ by $g=f_{i}'$ on each $[0,i].$ Now the function
$\varphi$ defined by 
$$\varphi(x)=\varphi_0+\int_0^{r}g(u)du$$ is a radial solution of (\ref{eqn k phi}).

\textit{Uniqueness:} assume that $\varphi_1$ and $\varphi_2$ are admissible radial solutions of (\ref{eqn k phi}): $\varphi_1(x)=f_1(r)$,
$\varphi_2(x)=f_2(r)$ where $f_1,f_2$ are functions $\R^+\rightarrow\R.$ For each real $R>0$, set
$$\varphi_{1,R}(x)=-\int_{r}^R{f_1}'(u)\ du\ \mbox{ and } \varphi_{2,R}(x)=-\int_{r}^R{f_2}'(u)\ du\ .$$
The functions $\varphi_{1,R}$ and $\varphi_{2,R}$ are both admissible solutions of the Dirichlet problem (\ref{radial Dirichlet problem i}) on $B_R.$ As
such, they must coincide on
$B_R,$ hence ${f_1}'={f_2}'$ on $[0,R],$ which implies the desired result.
\end{dem}
\subsection{Properties of the radial solutions}\label{parag 3}
The following lemma describes the monotonicity of a solution $\varphi$ of equation (\ref{eqn k phi}) depending on the sign of $h-1:$
\begin{Lemma}\label{variation phi} Let $h:\R^+\rightarrow\R$ and $\varphi:\H\rightarrow\R$ be as in Theorem \ref{theorem entire solution}, and let
$f:\R^+\rightarrow\R$ be such that $\varphi(x)=f[r(x)],$ $\forall x\in\H.$
\begin{enumerate}
\item[$(i)$] If $h\leq 1,$ then $f$ is non-increasing; in particular, if $\varphi_0=0$, the function $\varphi$ is non-positive.
\item[$(ii)$] If $h\geq 1,$ then $f$ is non-decreasing; in particular, if $\varphi_0=0$, the function $\varphi$ is non-negative.
\end{enumerate}
\end{Lemma}
\begin{dem}
Here, we need to calculate explicitely the expression of equation (\ref{eqn k phi}) in the radial case. Fix $x\in\H$ and take, with no loss of generality,
$$o=e_{n+1}=(0,\ldots,0,1),\
x=(\sinh r,0,\ldots,0,\cosh r)$$
with $r$, the hyperbolic distance between $o$ and $x.$ Consider the orthonormal basis of $T_x\H$
defined by:
$$\partial_r=\cosh r\ e_1+\sinh r\ e_{n+1},\mbox{ and }\partial_{\theta}=e_{\theta},\ \theta=2,\ldots,n,$$
and the vectors, tangent to $M=\graph_{\H}\phi$ at $e^{\varphi(x)}x$, induced by the embedding $x\in \H \to e^{\varphi(x)}x\in M$, given by: 
$$u_r=e^f(f'x+\partial_r),\ u_{\theta}=e^{f}\partial_{\theta},\ \theta=2,\ldots,n.$$
The future oriented unit normal to $M$ at $e^{\varphi(x)}x$ is the vector:
\begin{equation}\label{expr N}
N(r)=\frac{f'}{\sqrt{1-{f'}^2}}\ \partial_r+\frac{1}{\sqrt{1-{f'}^2}}\ x\ .
\end{equation}
Let $S$ be the curvature endomorphism of $M$ at $e^{\varphi(x)}x,$ with respect to the future unit normal $N(r).$ Using the formulas
$$D_{\partial_r}\partial_r=x,\ D_{\partial_\theta}\partial_r=\frac{1}{\tanh r}\partial_{\theta}$$
where $D$ denotes the canonical flat connection of $\R^{n,1},$ we readily get:
$$S(u_r)=dN(\partial_r)=\frac{e^{-f}}{\sqrt{1-{f'}^2}}\left(\frac{f''}{1-{f'}^2}+1\right)u_r\ ,$$
and, for $\theta=2,\ldots,n,$
$$S(u_{\theta})=dN(\partial_\theta)=\frac{e^{-f}}{\sqrt{1-{f'}^2}}\left(\frac{f'}{\tanh r}+1\right)u_\theta\ .$$
The principal curvatures of $M$ at $r>0$ are thus equal to:
$$\frac{e^{-f}}{\sqrt{1-{f'}^2}}\left(\frac{f''}{1-{f'}^2}+1\right)\mbox{ (simple)},\frac{e^{-f}}{\sqrt{1-{f'}^2}}\left(\frac{f'}{\tanh r}+1\right)\mbox{ (multiplicity
$n-1$)}.$$
Setting $s=s(r)$ for the hyperbolic distance from $o$ to $N(r),$ we infer from (\ref{expr N}):  
\begin{equation}\label{s fonction de f}
s(r)=r+\argth(f').
\end{equation}
In terms of the new radial unknown $s(r)$, for $r>0$, the principal curvatures reads
\begin{equation}\label{principal curvatures s}
\left(e^{-f}\cosh(r-s)s',e^{-f}\frac{\sinh s}{\sinh r},\ldots,e^{-f}\frac{\sinh s}{\sinh r}\right)\ ,
\end{equation}
and the equation $F_2(\varphi)=h$ reads
\begin{equation}\label{eqn s}
2s'\cosh(r-s)\sinh r\sinh s=nh^2\sinh^2 r-(n-2)\sinh^2 s.
\end{equation}
We now prove the first statement of the lemma. Since $f'=\tanh(s-r),$ we must prove: $s\leq r$ on $[0,+\infty).$ Suppose first $h<1.$ Since
$s(0)=0$ and $s'(0)=h(0)<1$ (from (\ref{eqn s})), there exists $r_0>0$ such that $s\leq r$ on $[0,r_0].$ Moreover, we get from (\ref{eqn s}):
$$s'\leq\frac{1}{2\cosh(r-s)}\left(n\frac{\sinh r}{\sinh{s}}-(n-2)\frac{\sinh s}{\sinh r}\right).$$
We observe that the function $s(r)=r$ is a solution of the ODE:
$$s'=\frac{1}{2\cosh(r-s)}\left(n\frac{\sinh r}{\sinh{s}}-(n-2)\frac{\sinh s}{\sinh r}\right)$$
on $[r_0,+\infty).$ So the comparison theorem for solutions of ordinary differential equations implies $s\leq r$ on $[r_0,+\infty).$
Suppose only $h\leq 1,$ fix $A>0$ and consider $h_{\delta}=h-\delta,$ where
$\delta$ is some small positive constant such that $h_{\delta}>0$ on $[0,A]$.
 Denoting by $\varphi_{\delta}$ and $s_{\delta}$ the corresponding solutions
of (\ref{eqn k phi}) and (\ref{eqn s}) on the ball of radius $A,$ the function
$s_{\delta}-r$ is non-positive; we now prove that  $s_{\delta}-r$ converges
uniformly to $s-r$ as $\delta$ tends to zero, which will yield the desired
result. Set $B_A$ for the ball of radius $A$ in $\H$ and 
$U=\{\psi\in\ C^{2,\alpha}(\overline{B_A}),\ \psi+\varphi\mbox{ is admissible in }\overline{B_A},\ \psi_{|\partial B_A}=0\}$; consider the auxiliary map:
$$\Phi: \psi \in U \to \Phi(\psi):=F_2(\psi+\varphi) \in C^{\alpha}(\overline{B_A})\ .$$
Since $\Phi(0)=h$ and since, classically \cite{gt} (recalling (\ref{expresslinear})), the linearized map $d\Phi(0)$ is an isomorphism from $\{\xi\in\
C^{2,\alpha}(\overline{B_A}),\ \xi_{|\partial B_A}=0\}$ to
$C^{\alpha}(\overline{B_A})$, the inverse function theorem implies: $\forall \varepsilon>0,
\exists \delta_0>0, \forall\delta\in(0,\delta_0)$, the solution $\psi_{\delta}\in U$ of
$F_2(\psi_{\delta}+\varphi)=h_{\delta}$ satisfies
$|\psi_{\delta}|_{2,\alpha}\leq\varepsilon.$ Since
$\varphi_{\delta}=\psi_{\delta}+\varphi-\psi_{\delta}(o),$ we obtain
$|\varphi_{\delta}-\varphi|_{2,\alpha}\leq 2\varepsilon,$ which implies the
convergence of $\varphi_\delta$ to $\varphi$ in $C^1$ and thus the uniform
convergence of $s_{\delta}$ to $s.$
\\The proof of statement $(ii)$ is analogous and thus omitted
\end{dem}

Our next lemma provides a simple necessary and sufficient condition for an entire radial solution to be bounded.

\begin{Lemma}\label{limite phi finie}
Let $h:\R^+\rightarrow\R$ and $\varphi:\H\rightarrow\R$ be as in Theorem \ref{theorem entire solution}.
\\$(i)$ Assume $h\leq 1,$ and $\lim_{r\rightarrow \infty}h= 1.$ Then
$$\lim_{r(x)\rightarrow+\infty}\varphi(x)>-\infty\mbox{ if and only if }\int_0^{+\infty}(1-h)dr\mbox{ converges.}$$
$(ii)$ Assume $h\geq 1,$ and $\lim_{r\rightarrow \infty}h= 1.$ Then
$$\lim_{r(x)\rightarrow+\infty}\varphi(x)<+\infty\mbox{ if and only if }\int_0^{+\infty}(h-1)dr\mbox{ converges.}$$
\end{Lemma}
\begin{dem}
Let us prove statement $(i)$, thus assuming $h\leq 1,$ with $\lim_{r\rightarrow \infty}h= 1.$ 
We stick to the notations used in the proof of Lemma \ref{variation phi}. From (\ref{s fonction de f}), we get at once:
\begin{equation}\label{phi fonction s}
\varphi(x)=\varphi_0-\int_0^{r(x)}\tanh(u-s(u))du\ .
\end{equation}
Statement $(i)$ amounts to prove that $\displaystyle{\int_{0}^{+\infty}\tanh(u-s(u))du}$ converges if and only if
so does $\displaystyle{\int_0^{+\infty}(1-h)dr}$. We split the proof of this fact into five steps.
\\
\\\textit{Step 1:} the solution $s$ of (\ref{eqn s}) is an increasing function. 

Let us consider in the $(r,s)$ plane  the curve ${\cal C}$ with equation:
$$nh^2\sinh^2r= (n-2)\sinh^2s,\ r,s\geq 0.$$
The slope of its tangent at $(0,0)$ is $\sqrt{\frac{n}{n-2}}h(0).$ Since the solution $s$ satisfies $s(0)=0$ and $s'(0)=h(0),$ we infer that the graph of
$s$ stays under the curve ${\cal C}$ near 0. Noting that the following vector field, associated to the differential equation (\ref{eqn s}):
$$(r,s)\mapsto(2\cosh(r-s)\sinh r\sinh s,nh^2\sinh^2r-(n-2)\sinh^2s)\ ,$$
is horizontal on ${\cal C}$, and that the height $s$ of the curve ${\cal C}$ is increasing with $r$, we conclude that the solution $s$ of (\ref{eqn s})
remains trapped below
${\cal C}.$ In other words
$nh^2\sinh^2r\geq (n-2)\sinh^2s$ for all $r,$ and (\ref{eqn s}) implies: $s'\geq 0.$
\\ 
\\\textit{Step 2:}   $r-s$ has a limit at $+\infty.$ 

By contradiction, assume $\liminf(r-s)<\limsup(r-s)=\delta.$ Thus there exists a sequence $r_k\rightarrow+\infty$ such that
$r_k-s(r_k)\rightarrow\delta$ and $s'(r_k)=1.$ Denoting $s(r_k)$ by $s_k,$ we get from equation (\ref{eqn s}):
\begin{equation}\label{eqn sk}
1=\frac{1}{2\cosh(r_k-s_k)}\left[nh^2(r_k)\frac{\sinh r_k}{\sinh s_k}-(n-2)\frac{\sinh s_k}{\sinh r_k}\right].
\end{equation}
We distinguish two cases :
\\\textit{First case: $\delta<+\infty.$} We then have $s_k\rightarrow+\infty,$ $\frac{\sinh r_k}{\sinh s_k}\sim e^{r_k-s_k}\sim e^{\delta}$ and  $\frac{\sinh s_k}{\sinh r_k}\sim e^{s_k-r_k}\sim e^{-\delta}$ as $k$ tends to
infinity (here and below, the equivalence $\sim$ between two quantities means that their quotient has limit 1).
So (\ref{eqn sk}) yields
$$1=\frac{1}{2\cosh\delta}\left[ne^\delta-(n-2)e^{-\delta}\right].$$
Using $e^\delta\geq e^{-\delta}$ we get $1\geq\frac{e^{\delta}}{\cosh\delta},$ which is absurd.
\\
\\\textit{Second case $\delta=+\infty.$} First assuming that $s_k$ is not bounded, and since $s$ is an increasing function (Step 1), we have :
$s_k\rightarrow+\infty,$ $\frac{\sinh r_k}{\sinh s_k}\sim e^{r_k-s_k}\rightarrow+\infty$ and  $\frac{\sinh s_k}{\sinh r_k}\sim e^{s_k-r_k}\rightarrow 0$ as
$k$ tends to infinity. Equation (\ref{eqn sk}) yields 
$$1\sim\frac{n}{2\cosh(r_k-s_k)}e^{r_k-s_k},$$
which is absurd since $\cosh(r_k-s_k)\sim\frac{e^{r_k-s_k}}{2}.$ If we now assume $s_k$ bounded, since $s$ is an increasing function with $s'(0)>0$,
we get that $s_k$ converges to $l>0,$ and, since $\frac{\sinh s_k}{\sinh r_k}\rightarrow 0,$ we obtain from (\ref{eqn sk}):
$$1\sim\frac{n}{2\cosh(r_k-s_k)}\frac{\sinh r_k}{\sinh l},$$
with $\sinh r_k\sim\frac{e^{r_k}}{2},$ $\cosh(r_k-s_k)\sim\frac{e^{r_k-s_k}}{2}\sim \frac{e^{-l}}{2}e^{r_k};$ so $1=\frac{n}{2}\frac{e^l}{\sinh l},$ which
is absurd. 
\\
\\\textit{Step 3: } $r-s$ tends to 0 at infinity.

Having proved that $r-s$ converges, let us set $\delta=\lim_{r\rightarrow+\infty}r-s$ and prove by contradiction that $\delta=0.$ There are two cases :
\\
\\\textit{First case : $0<\delta<+\infty.$} We get $s\rightarrow+\infty,$ hence $\frac{\sinh r}{\sinh s}\sim e^{r-s}\sim e^{\delta},$ $\frac{\sinh s}{\sinh
r}\sim e^{s-r}\sim e^{-\delta}$ as $r$ tends to infinity, and thus, from (\ref{eqn s}):
$$s'\rightarrow \frac{1}{2\cosh\delta}\left[ne^{\delta}-(n-2)e^{-\delta}\right].$$
The latter expression is larger than 1, which contradicts $r\geq s.$
\\
\\\textit{Second case : $\delta=+\infty.$} We first note that $\frac{\sinh s}{\sinh r}\rightarrow 0$ (if $s$ is bounded this is trivial; if $s$ is not
bounded,
$s\rightarrow+\infty$ since $s$ is increasing, and we have $\frac{\sinh s}{\sinh r}\sim e^{s-r}\rightarrow 0$ since $r-s\rightarrow+\infty$). Moreover we
have $\liminf nh^2\frac{\sinh r}{\sinh s}\geq n$ since $r\geq s.$ We thus infer from equation (\ref{eqn s}):
$$s'\sim\frac{n}{2\cosh(r-s)}\frac{\sinh r}{\sinh s}.$$
Assuming $s\rightarrow+\infty,$ we get $\frac{\sinh r}{\sinh s}\sim e^{r-s}$ and $\cosh(r-s)\sim\frac{e^{r-s}}{2},$ hence $s'\sim n,$ which is
impossible since $s\leq r.$
\\Finally, assuming $s$ bounded yields $s\rightarrow l>0;$ since $r-s\rightarrow+\infty,$ we infer $\cosh(r-s)\sim \frac {e^{r-s}}{2}$ and $\sinh r\sim
\frac{e^r}{2},$ hence from (\ref{eqn s}), $e^{-s}s'\sim\frac{n}{2}\frac{1}{\sinh l}$ and thus $s'\sim \frac{n}{2}\frac{e^l}{\sinh l},$ which contradicts
the boundedness assumption on $s$.
 \\
 \\\textit{Step 4:} $\lim_{r(x)\rightarrow+\infty}\varphi(x)>-\infty$ if and only if $\varepsilon(r):=r-s$ is integrable on $[0,+\infty).$ 
 
This is straightforward from (\ref{phi fonction s}) combined with $\tanh(u-s(u))\sim\varepsilon(u)$ which holds as $u\to +\infty$ due to Step 3.
\\
\\\textit{Step 5:} $\varepsilon$ is integrable on $[0,+\infty)$ if and only if $\beta:=1-h^2$ is integrable on $[0,+\infty).$

First observation: $\lim_{r\to \infty}s'=1$. Indeed, at infinity, we have $r-s\rightarrow 0$, so $s\rightarrow+\infty,$ hence: 
$$\frac{\sinh r}{\sinh s}\sim e^{r-s}\sim 1,\ \frac{\sinh s}{\sinh r}\sim e^{s-r}\sim 1,$$
and (\ref{eqn s}) yields $s'\rightarrow 1.$
\\

Using Step 3, the assumptions on $h$ and the preceding observation, we get
$$\varepsilon(r)\rightarrow 0,\ \beta(r)\rightarrow 0,\mbox{ and }\varepsilon'(r)=1-s'(r)\rightarrow 0$$
as $r$ tends to infinity. Plugging the definitions of $\varepsilon$ and $\beta$ in (\ref{eqn s}) and using the expansions 
$$\cosh \varepsilon=1+o(\varepsilon),\ \sinh(r-\varepsilon)=\sinh r\ (1-\varepsilon+o(\varepsilon)),$$
yields
\begin{equation}
(n-1)\varepsilon+\varepsilon'+o(\varepsilon)=\frac{n}{2}\beta.
\end{equation}
Fixing a real  $\delta>0,$ there readily exists $r_{\delta}>0$ such that, for all $r\geq r_{\delta},$ 
\begin{equation}\label{ineq epsilon 1}
\varepsilon'+(n-1-\delta)\varepsilon\leq\frac{n}{2}\beta\ ,
\end{equation}
and
\begin{equation}\label{ineq epsilon 2}
\varepsilon'+(n-1+\delta)\varepsilon\geq\frac{n}{2}\beta\ .
\end{equation}
Integrating (\ref{ineq epsilon 1}), we get, for $r\geq r_{\delta},$
$$\varepsilon(r)\leq e^{-(n-1-\delta)r}\left[C(r_{\delta})+\frac{n}{2}\int_{r_{\delta}}^r\beta(u)e^{(n-1-\delta)u}du\right].$$
Integrating again and using Fubini Theorem yields, with $\delta$ such that $n-1-\delta>0,$
\begin{eqnarray*}
\int_{r_{\delta}}^{+\infty}\varepsilon(r)dr&\leq& C'(r_{\delta})+\frac{n}{2}\int_{r_{\delta}}^{+\infty}\beta(u)e^{(n-1-\delta)u}\left(\int_{u}^{+\infty}e^{-(n-1-\delta)r}dr\right)du,\\
&\leq & C'(r_{\delta})+\frac{n}{2(n-1-\delta)}\int_{r_{\delta}}^{+\infty}\beta(u) du.
\end{eqnarray*}
We conclude that $\varepsilon$ is integrable provided $\beta=1-h^2$ is integrable.
\\Analogously, using (\ref{ineq epsilon 2}), we get 
$$\varepsilon(r)\geq e^{-(n-1+\delta)r}\left[C(r_{\delta})+\frac{n}{2}\int_{r_{\delta}}^r\beta(u)e^{(n-1+\delta)u}du\right],$$
and
\begin{eqnarray*}
\int_{r_{\delta}}^{+\infty}\varepsilon(r)dr&\geq& C'(r_{\delta})+\frac{n}{2}\int_{r_{\delta}}^{+\infty}\beta(u)e^{(n-1+\delta)u}\left(\int_{u}^{+\infty}e^{-(n-1+\delta)r}dr\right)du,\\
&\geq & C'(r_{\delta})+\frac{n}{2(n-1+\delta)}\int_{r_{\delta}}^{+\infty}\beta(u) du.
\end{eqnarray*}
Taking $\delta>0$ arbitrary, we find that $\beta$ is integrable if $\varepsilon$ is integrable.
\\The proof of statement $(ii)$ is analogous and thus omitted
\end{dem}
\subsection{Construction of appropriate radial barriers}\label{parag 4}
\begin{Lemma}
Let $h:\H\rightarrow\R$ be a positive and continuous function on the hyperbolic space such that
$$\lim_{r(x)\rightarrow+\infty} h(x)=1$$
and such that  the functions $h^-$ and $h^+$ defined on $\R^+$ by
$$h^-(r)=\sup_{r(x)=r} h(x)\mbox{ and }h^+(r)=\inf_{r(x)=r} h(x)$$
satisfy 
$$ \int_0^{+\infty}(h^--1)_+dr<+\infty,\ \int_0^{+\infty}(1-h^+)_+dr<+\infty,$$  
where $(h^--1)_+$ (resp. $(1-h^+)_+$) means the positive part of $h^--1$ (resp. $1-h^+$). Then there exist $\varphi^-,\varphi^+\in C^{\infty}(\H)$, with
strictly convex spacelike graphs, satisfying:
$$F_2(\varphi^-)\geq h,\ F_2(\varphi^+)\leq h\mbox{ and }\lim_{r\rightarrow+\infty}\varphi^{\pm}=0.$$
\end{Lemma}
\begin{dem}
First, considering $1+(h^--1)_+$ instead of $h^-$ and $1-(1-h^+)_+$ instead of $h^+,$ we may suppose without loss of generality that $h^-$ and $h^+$ are two
continuous functions such that : $\forall x\in\H,$ with $r=r(x)$,
\begin{equation}\label{cond hpm 1}
h^-(r)\geq h(x)\geq h^+(r)>0,
\end{equation}
\begin{equation}\label{cond hpm 2}
h^-\geq 1\geq h^+,\ \lim_{r\rightarrow+\infty}h^-(r)=\lim_{r\rightarrow+\infty}h^+(r)=1,
\end{equation}
and 
\begin{equation}\label{cond hpm 3}
\int_0^{+\infty}(h^--1)dr<+\infty,\ \int_0^{+\infty}(1-h^+)dr<+\infty.
\end{equation} 
If we now consider 
$$h^-+\frac{\varepsilon_0}{r^2}\mbox{ if }r\geq 1,\ h^-+\varepsilon_0\mbox{ if }r\leq 1$$
instead of $h^-,$ and 
$$h^+-\frac{\varepsilon_0}{r^2}\mbox{ if }r\geq 1,\ h^+-\varepsilon_0\mbox{ if }r\leq 1$$ instead of $h^+,$ where $\varepsilon_0$ is chosen sufficiently small such that $\inf
h^+>\varepsilon_0,$
we may moreover assume the following:
$$h^-\geq\max(1,h)+\frac{\varepsilon_0}{r^2}\mbox{ and } h^+\leq\min(1,h)-\frac{\varepsilon_0}{r^2}\mbox{ if }r\geq 1.$$
We now prove that we can approximate $h^\pm$ by smooth functions $g^\pm$ such that 
\begin{equation}\label{estimate approximation}
|h^\pm-g^\pm|\leq\min\left(\frac{\varepsilon_0}{r^2},\varepsilon_0\right).
\end{equation}
 For each $i\in\N,$ let us denote by $g^-_i$ a smooth function on $[0,i+1]$ such that $|h^--g^-_i|\leq\frac{\varepsilon_0}{(i+1)^2}$ on $[0,i+1].$ Let $\theta\in
C^{\infty}_c(\R)$
such that $0\leq\theta\leq 1,$ $\theta(x)=1$ if $|x|\leq\frac{1}{4}$ and $\theta(x)=0$ if $|x|\geq\frac{3}{4}.$ We define $g^-$ on
$[i,i+1]$ by 
$$g^-=\theta_ig^-_i+(1-\theta_i)g^-_{i+1},$$
where $\theta_i=\theta(.-i).$ By construction, we have $g^-=g^-_i$ on a neighborhood of $i.$ The function $g^-$ is thus smooth on $[0,+\infty),$ and satisfies on $[i,i+1]$ :
$$|g^--h^-|\leq \theta_i|g^-_i-h^-|+(1-\theta_i) |g^-_{i+1}-h^-|\leq\frac{\varepsilon_0}{(i+1)^2},$$
which implies the estimate (\ref{estimate approximation}). We may thus assume that (\ref{cond hpm 1}), (\ref{cond hpm 2}) and (\ref{cond hpm 3}) hold, where
$h^{\pm}$ are two smooth functions on $\R^+.$ Considering $\theta\sup_{\R^+}h^-+(1-\theta)h^-$ instead of $h^-,$ and $\theta\inf_{\R^+}h^++(1-\theta)h^+$
instead of $h^+,$ we may also assume that the functions ${h^\pm}$ are constant on some neighborhood of 0. Let $\varphi^-$ and $\varphi^+$ be smooth radial
functions given by Theorem \ref{theorem entire solution} (with some arbitrary initial condition $\varphi_0$) such that $F_2(\varphi^\pm)=h^\pm$. From Lemma
\ref{limite phi finie}, subtracting constants if necessary, we obtain $\lim_{r\rightarrow+\infty}\varphi^\pm(r)=0$ 
\end{dem}
\end{document}